# The focused information criterion for stochastic model selection problems using $M$-estimators

BY S. C. PANDHARE AND T. V. RAMANATHAN

*Department of Statistics and Center of Advanced Studies, Savitribai Phule Pune University, Pune-411007, India*

pandhare.05@gmail.com ram@unipune.ac.in

SUMMARY

Claeskens and Hjort (2003) constructed the focused information criterion FIC and developed frequentist model averaging methods using maximum likelihood estimators assuming the observations to be independent and identically distributed. Towards the immediate extensions and generalizations of these results, the present article is aimed at providing the focused model selection and model averaging methods using general maximum likelihood type estimators, popularly known as $M$-estimators. The necessary asymptotic theory is derived in a setup of stationary and strong mixing stochastic processes employing von Mises functional calculus of empirical processes and Le Cam's contiguity lemmas. We illustrate the proposed focused stochastic modeling methods using three well-known spacial cases of $M$-estimators, namely, conditional maximum likelihood estimators, conditional least square estimators and estimators based on method of moments. For the sake of simulation exercises, we consider two simple applications of FIC. The first application discusses the simultaneous selection of order of autoregression and symmetry of innovations in asymmetric Laplace autoregressive models. The second application demonstrates the FIC based choice between general scale-shape Gamma density and exponential density with shape being unity. We observe that in terms of the correct selections, FIC outperforms classical Akaike's information criterion AIC and performs at par with Bayesian information criterion BIC.

## 1. INTRODUCTION

The search for the best possible and most explanatory model for the given data is an integral part of any data analysis and decision making process. The context in which the models are to be selected and their role in explaining the associations and interrelations among the variables in the given data are usually determined by certain functions of parameters of interest. Hence, the criterion that is to be employed to select a model from the given set of the models must take into consideration this functional form of parameters of interest. The classical model selection criteria like Akaike's information criterion AIC and Bayesian information criterion BIC provide a single, most suited model for the data based on the global fit of the candidate models. The focused information criterion or FIC of Claeskens and Hjort (2003) on the other hand is based on asymptotic mean squared error of estimators of user-defined focus functions chosen according to the context and purpose of models.

Ever since its inception, this idea of focused model selection of Claeskens and Hjort and inference has received considerable attention and appreciations. Hansen (2005) describes FIC to be " an intriguing challenger to existing model selection methods and deserves attention and scrutiny". Focused model selection methods have been developed under a variety of data generating mechanisms such as regression models (Claeskens and Hjort, 2008), autoregressive time series models (Claeskens et al., 2007, Rohan and Ramanathan, 2011), semiparametric generalized lienar models (Claeskens and Caroll, 2007) etc. Claeskens (2016) provides an extensive review of the recent advances in this domain of statistical modeling and inference.

The focused model selection and model averaging methods of Claeskens and Hjort (2003) rely on the likelihood based estimation and assume the independence of observations. Although use of efficient likelihood based estimators is widespread, these estimators are highly sensitive to even moderate deviations in the hypothetical, "true" probability measure assumed to have been fixed a priori. The local perturbation of the form $\theta_n = \theta_0 + n^{-1/2}\delta$ in the sequence of fixed probability measures $P_{n,\theta_0}$ is crucial to the focused model selection and model averaging theory of Claeskens and Hjort (2003). It is possible to establish the asymptotic normality of sequence of likelihood based estimators even in presence of such local deviations. Le Cam's local asymptotic normality lemmas ascertain that these estimators do attain asymptotic Cramer-Rao lower bound with the variance of asymptotic normal distribution of scaled estimators being the inverse of the Fisher information matrix. However, the estimators are no longer asymptotically unbiased owing to the local misspecification. The focused model selection typically involves comparisons of the mean squared errors of the asymptotically normal estimators of user-defined focus functions. It is therefore essential to take into consideration the trade-off between efficiency and bias while selecting models under the local misspecification. Developing model selection and inference methods by means of perhaps non-efficient estimators but which are in some sense robust to such distributional misspecification is worth the investigations. With this view, we propose the focused model selection using general maximum likelihood-type estimators or simply $M$-estimators of Huber (1981) nesting the important special case of efficient likelihood estimators.

Given a random sample $X_1, ..., X_n$ of size $n$, $M$-estimators are viewed as the optima of certain criterion function, say, $\sum_{i=1}^{n} M(X_i, \theta)$. Under moderate differentiability assumptions on the criterion function, the $M$-estimators may also be viewed at as consistent roots of the estimating equations $\sum_{i=1}^{n} \psi(X_i, \theta) = 0$ where $\psi(., \theta)$ denotes the gradient of $M(., \theta)$. It is clear that in principle, the criterion function may not be specified in terms the distributional attributes such as probability density functions. Hence, these estimators are expected to safeguard against the possible distributional incorrectness. The details on the robustness of $M$-estimators are referred to Hampel (1974) and Hampel et al. (2011).

To the best of our knowledge, there are only two articles that propose focused model selection and model averaging using non-likelihood estimators, both of which in their own way form the special cases of the approach discussed here. DiTraglia (2016) develops the FIC using generalized method of moments for the focused selection of moments from a set of potentially misspecified moment conditions. The author derives the asymptotic theory under the misspecification of the form $E\left(\psi\left(X, \theta\right)\right) = n^{-1/2}\eta$ assuming numerous stringent regularity conditions. The main limitation of these discussions lies in the fact that moments seldom assist in explaining the data generating mechanisms. Statistical model selection in a broad sense is essentially a selection of an appropriate probability measure that may explain the data generating mechanism to a desired accuracy. Hence, the methods of DiTraglia may not be applicable for general model selection problems. On the other hand, like Claeskens and Hjort (2003), we assume the parametric misspecification in probability measures which implicitly also suggests that the estimating functions (the moment conditions of DiTraglia, 2016) are locally misspecified. The converse implication does not hold true in general. Moreover, our local misspecification setup ensures the contiguity of misspecified probability measures with the fixed probability measure allowing one to utilize

Le Cam's lemmas leading to the transparent proofs of asymptotic theory. The asymptotic theory of empirical process and functional calculus enable us to derive the local asymptotic theory of $M$-estimators with much weaker, fewer assumptions than those of DiTraglia.

The key result as presented in our Theorem 2 closely follows the proof of local asymptotic normality of $M$-estimators for independent, identically distributed observations discussed by Rieder (2012) Chapter 4, Lemma 4.2.18. In particular, the asymptotic normality as claimed in our Theorem 2 parallels the equation (66) of Rieder (2012), Chapter 4. We resort to Le Cam's contiguity lemmas extending the results of Rieder to dependent setup with applications to model selection problems.

The other instance of using non-likelihood estimators for focused model selection is provided by Lohmeyer et al. (2018) who construct FIC using least square estimators for selection problems for vector autoregressive models. Least square estimators are indeed $M$-estimators with the squared residual function being the criterion function to be minimized. Although we present the results for univariate stochastic process, it should be possible to extend the functional asymptotics and consequently the parametric inference to the case of vector valued stochastic processes. Hence, it is possible to view the theory of Lohmeyer et al. (2018) as special case of our setup.

In a nutshell, we aim for the generalizations of the focused model selection theory to the stochastic dependence setup, developing the FIC using general $M$-estimators and evaluating empirically the effects of the bias-variance trade-off as exhibited by these estimators on the accuracy of model selection. We illustrate our method using estimators based on conditional maximum likelihood, conditional least squares and method of moments. However, it is possible to obtain the FIC expressions using any other, perhaps more robust $M$-estimators by virtue of our general results.

## 2. Preliminaries

### 2·1. *Notations and assumptions*

Let $\{X_n, n \in \mathbb{Z}^+\}$ denote a stochastic process on the probability space $(\Omega, \mathscr{F}, P_\theta)$ where $\mathbb{Z}^+$ denote the set of non-negative integers. Throughout the subsequent discussions, $\theta \in \Theta \subseteq \mathbb{R}^k$ denotes the labeling parameter of the underlying probability models. Without loss of generality, we assume that $\{X_n\}$ is a continuous state space stochastic process. Further, let $\mathscr{G}$ and $\mathscr{H}$ denote sub-$\sigma$-fields of $\mathscr{F}$. Define an $\alpha$-mixing coefficient as

$$\alpha\left(\mathscr{G}, \mathscr{H}\right) = \sup_{A \in \mathscr{G}, B \in \mathscr{H}} \left|P\left(A \cap B\right) - P(A)P(B)\right|.$$

Denote a sequence $\{\alpha_n, n \in \mathbb{Z}^+\}$ as

$$\alpha_n = \alpha\left(\sigma\left(X_0\right), \sigma\left(X_n, X_{n+1}, ...\right)\right).$$

Let $E(X)$ denote the expectation of a random variable $X$ under $P_\theta$. The following are the fundamental assumptions that the present technical discussions rely on.

*Assumption* 1 (*A1*). $\{X_n, n \in \mathbb{Z}^+\}$ is a strictly stationary stochastic process.

*Assumption* 2 (*A2*). $E(X_0) = 0$ without loss of generality and $E|X_0|^b < \infty$ for some $b > 2$.

*Assumption* 3 (*A3*). $X_0$ is bounded with $P_\theta$-probability 1.

*Assumption* 4 (*A4*). $\sum_{n=0}^{\infty} n^2 \alpha_n^{1-2/b} < \infty$ for some $b > 2$.

Let $\theta_0$ be a fixed interior point of $\Theta$. Suppose $F_{\theta_0}$ is the proposed finite dimensional distribution of a random sample $X_1, ...., X_n$. Let $F_n$ denote the empirical distribution function of this sample.

Set $\mathcal{L}^2(F_\theta) = \{h : \Omega \times \Theta \to \mathbb{R}^k \mid E||h(X,\theta)h(X,\theta)^T|| < \infty\}$ where $||\ ||$ denotes the matrix norm. Let $D[-\infty,\infty]$ denote a collection of right continuous functions on $[-\infty,\infty]$ with finite left hand limits equipped with the uniform topology. Let $T : D[-\infty,\infty] \to \mathbb{R}^k$, denote a statistical functional. We assume that $T$ is a Fisher-consistent functional, meaning that $T(F_\theta) = \theta_0$. Let $0_m$ denote a vector of zeros of length $m$.

Now, consider a measurable function $\psi : \Omega \times \Theta \to \mathbb{R}^k, \Omega \subseteq \mathbb{R}$ such that $\psi(x,.)$ is almost surely continuous and $\psi(.,\theta) \in \mathscr{L}^2(F_\theta)$ for every $\theta \in \Theta$. Following Huber (1981), an $M$-functional is defined as a root $T(F_\theta) = \theta_0$ of the equation

$$\int \psi(x,\theta)F_\theta(dx) = E\left(\psi(X_0,\theta)\right) = 0_k. \tag{1}$$

Corresponding to $M$-functional, an $M$-estimator $T(F_n) = \hat{\theta}$ is defined as a root of the equation

$$n^{-1}\sum_{t=1}^{n} \psi\left(X_t,\theta\right) = 0_k. \tag{2}$$

Set $\phi(\theta) = E\left(\psi(X_0,\theta)\right)$. In order to derive the asymptotics of $M$-estimators, the following is a crucial regularity condition.

*Assumption* 5 (*A5*). The function $\phi$ is homeomorphic at $\theta$ and has invertible, continuous, bounded derivative (Jacobian matrix) $d\phi(\theta)$ at $\theta$.

Let $\mathscr{F}_{n-1} = \sigma\left(X_{n-1},...,X_0\right)$ denote the $\sigma$-field generated by $X_{n-1},...,X_0$. Let $f_\theta\left(x_0\right)$ and $f_\theta\left(x_n|\mathscr{F}_{n-1}\right)$ denote respectively the initial and transition density functions. Let $\{\theta_n, n \in \mathbb{Z}^+\}$ denote a sequence in $\Theta$. Consider the sequence of likelihood ratios

$$\Lambda_n\left(\theta_0,\theta_n\right) = \frac{f_{\theta_n}\left(x_0\right)}{f_{\theta_0}\left(x_0\right)}\prod_{t=1}^{n}\frac{f_{\theta_n}\left(x_t|\mathscr{F}_{t-1}\right)}{f_{\theta_0}\left(x_t|\mathscr{F}_{t-1}\right)}. \tag{3}$$

Under the local misspecification setting, $\theta_n = \theta_0 + n^{-1/2}\delta$ for some $\delta \in \mathbb{R}^k$.

*Assumption* 6 (*A6*). There exists a measurable function $S : \Omega \times \Theta \to \mathbb{R}^k$ such that

$$n^{1/2}\left(\frac{f_{\theta_n}^{1/2}\left(x_n|\mathscr{F}_{n-1}\right)}{f_{\theta_0}^{1/2}\left(x_n|\mathscr{F}_{n-1}\right)} - 1\right) - \frac{1}{2}\delta^{\mathrm{T}}S(X_n,\theta_0) \to 0 \text{ in quadratic mean.}$$

### 2·2. *Discussions of assumptions*

The local misspecification is required to prevent the bias of estimates to explode when the estimation is carried out under possibly wrong models. Besides, the local misspecification assumption that we work with is deeply intertwined with Le Cam's contiguity of probability measures. We refer to Rousass (1972) and Le Cam and Yang (2012) for the details on contiguity and its statistical applications. In simple words, the local misspecification assumes that the "true" model lies in a vicinity of the (possibly wrong) model fixed a priori. The inference should then be carried in this neighbourhood. The contiguity or the nearness of these misspecified models to the fixed models makes such inference possible while safeguarding the inference procedures against the consequences of working with fixed, wrong models.

Assumptions A1 and A4 together imply that the stochastic process is strong $(\alpha)$-mixing. We refer to Bradley (2005) for the properties of various mixing conditions. The assumptions A1-A4 imply that the central limit theorem holds true even in the dependent setup, meaning that the summand $n^{-1/2}\sum_{t=1}^{n} X_t$ converges in distribution under $P_\theta$ to normal random variable with mean 0 and finite variance $\sigma^2$, cf. Deo (1973), Lemma 2. The said central limit theorem along with the tightness lemma of Goldie and Greenwood (1986) and metric entropy central limit theorem of Ossiander (1987) is instrumental in establishing the weak convergence of empirical process $n^{1/2}\left(F_n - F_{\theta_0}\right)$ to the tight Brownian motion. Since $T$ is assumed to be Fisher consistent, the asymptotics of this empirical process can be connected to those of sequence of estimators of $\theta$ by means of the functional delta method of Fernholz (1983) provided $T$ is a Hadamard differentiable functional. These functional asymptotics are discussed in Appendix.

Assumption A6 implies that the underlying density is $\mathscr{L}^2$-differentiable or differentiable in quadratic mean (Rieder, 2012, Chapter 2, Definition 2.3.6.). The quadratic mean derivative $S(X_n, \theta_0)$ in A6 coincides with the conditional score vector $\frac{\partial \log f_\theta(x_n|\mathscr{F}_{n-1})}{\partial \theta_0}$. The assumption of $\mathcal{L}^2$-differentiability requires the existence of only the first order partial derivatives of conditional density except on a set of Lebesgue measure zero. Since no assumptions on the existence of higher order derivatives of the density are needed, this can be seen to be much weaker assumption than the standard Cramer-regularity conditions that are often assumed to derive the asymptotics of the likelihood-based estimators. See, for instance, Basawa and Rao (1980). We deviate from this routine and present the results in a generalized setup of functional asymptotics. With the estimating function $\psi\left(X_t, \theta\right)$ in (1) and (2) above set equal to the score $S\left(X_t, \theta\right)$, A6 implies A5 with $d\phi(\theta) = E\left(S(X_0, \theta) S(X_0, \theta)^{\mathrm{T}}\right) = J(\theta)$, the Fisher information matrix. Hence, this route of functional calculus and $M$-estimators that we undertake also reestablishes the local asymptotic normality of maximum likelihood estimators with much weaker differentiability assumptions on likelihood than those of Hjort and Claeskens (2003).

## 3. Asymptotic Theory

LEMMA 1. *Under the sequence of fixed probability measures $P_{n,\theta_0}$ and assumptions A1-A5, an $M$-functional $T$ defined by (1) is asymptotically linear with the influence curve $-d\phi\left(\theta_0\right)^{-1} \psi\left(x, \theta_0\right)$.*

*Proof.* Since the process $\{X_n\}$ is stationary and strong mixing, the process $\{F_\theta(X_n)\}$ is again stationary and strong mixing on $[0,1]$ for every $\theta$. Thus, on the lines of Fernholz (1983), Subsections 4.1 and 5.1, it can be seen that the $M$-functional (1) induces implicitly another functional $\tau : D[0,1] \to \mathbb{R}^k$ such that for any continuous and strictly increasing distribution function $F^* \in D[0,1]$, $\tau(F^*) = T(F^* \circ F)$. Hence, $\tau(F^*) = \theta$ is a root of $\rho\left(F^*, \theta\right) = 0_k$ where an $i^{\text{th}}$ component of vector $\rho$ is given as

$$\rho_i\left(F^*, \theta\right) = \int_0^1 \psi_i\left(F_\theta{}^{-1}\left(F^{*-1}(x)\right), \theta\right) dx, \ i = 1, 2, ..., k.$$

Assumption A5 implies that the conditions of Proposition 7.1.1 of Fernholz (1983) are satisfied implying the Hadamard differentiability of $\rho$ at $(F^*(X_0), \theta_0)$. Now invoke the implicit function theorem for Hadamard differentiable functionals (Fernholz 1983, Theorem 6.2.1) to conclude that $M$-functional is Hadamard differentiable at $F_\theta$. Identifying $T(F_n)$ as $\hat{\theta}$ and $T(F_{\theta_0})$ as $\theta_0$, the asymptotic linearity follows due to Theorem A1 in Appendix. The derivation of influence curve is referred to Huber (1981), Hampel (1974). □

Next, we present the key result on asymptotic expansion of likelihood ratio (3). The proof is provided in the supplementary material.

THEOREM 1. *Under A1-A4 and A6,*

$$\log \Lambda_n(\theta_n, \theta_0) = n^{-1/2}\delta^{\mathrm{T}} \sum_{t=1}^{n} S(X_t, \theta_0) - \frac{1}{2}\delta^{\mathrm{T}} J(\theta_0)\delta + o_{P_{\theta_0}}(1).$$

Finally, the following result provides the local asymptotic normality of $M$-estimators useful to derive the FIC expressions.

THEOREM 2. *Under the assumptions A1-A6,*

$$n^{1/2}(\hat{\theta} - \theta_0) \to N_k\left(-b(\delta), D\Sigma D^{\mathrm{T}}\right) \text{ in distribution under local misspecification } P_{n,\theta_0+n^{-1/2}\delta}. \tag{4}$$

*Here,* $b(\delta) = DE\delta$, *where,* $D = d\phi(\theta_0)^{-1}$, $\Sigma = E\left(\psi\left(X_0, \theta_0\right)\psi\left(X_0, \theta\right)^{\mathrm{T}}\right)$ *and* $E = E\left(\psi\left(X_0, \theta_0\right) S\left(X_0, \theta\right)^{\mathrm{T}}\right)$.

*Proof.* Lemma 1 implies that

$$n^{1/2}(\hat{\theta} - \theta_0) = -n^{-1/2} D \sum_{t=1}^{n} \psi\left(X_t, \theta_0\right) + o_{P_{\theta_0}}(1).$$

In view of Theorem A1 in Appendix, the above relation implies

$$n^{1/2}\left(\hat{\theta} - \theta_0\right) \to N_k\left(0, D\Sigma D^{\mathrm{T}}\right) \text{ in distribution.} \tag{5}$$

Since stationary, strong mixing process is ergodic (Bradley, 2005, Subsection 2.5), Assumption A6 implies that the sequence of score vectors $\{S(X_n, \theta)\}$ is a martingale difference sequence, cf. Lemma 2 in the supplementary material. The martingale central limit theorem (Billingsley, 1961b) coupled with Cramer-Wold device (van der Vaart, 2000, page 16) now implies that for any $\delta \in \mathbb{R}^k$,

$$n^{-1/2}\delta^{\mathrm{T}} \sum_{t=1}^{n} S(X_t, \theta_0) \to N_1\left(0, \delta^{\mathrm{T}} J(\theta_0)\delta\right) \text{ in distribution.}$$

Slutsky's lemma (van der Vaart, 2000, Lemma 2.8)) in view of Theorem 1 now shows that under $P_{n,\theta_0}$,

$$\log \Lambda_n(\theta_n, \theta_0) \to N_1\left(-\frac{1}{2}\delta^{\mathrm{T}} J(\theta_0)\delta, \delta^{\mathrm{T}} J(\theta_0)\delta\right) \text{ in distribution.} \tag{6}$$

Applying Cramer-Wold device to (5) and (6) now yields

$$\begin{pmatrix} n^{1/2}(\hat{\theta} - \theta_0) \\ \log\Lambda_n\left(\theta_n, \theta_0\right) \end{pmatrix} \to N_{k+1}\left(\begin{pmatrix} 0_{k+1} \\ -\frac{1}{2}\delta^{\mathrm{T}} J(\theta_0)\delta \end{pmatrix}, \begin{pmatrix} D\Sigma D^{\mathrm{T}} & -b(\delta) \\ -b(\delta)^{\mathrm{T}} & \delta^{\mathrm{T}} J(\theta_0)\delta \end{pmatrix}\right) \text{ in distribution.} \tag{7}$$

Le Cam's first lemma (van der Vaart, 2000, Lemma 6.4) in view of (7) implies the mutual contiguity of $P_{n,\theta_0}$ and $P_{n,\theta_0+n^{-1/2}\delta}$. Applying Le Cam's third lemma (van der Vaart, 2000, Lemma 6.6) shows that under $P_{n,\theta_0+n^{-1/2}\delta}$,

$$\begin{pmatrix} n^{1/2}\left(\hat{\theta} - \theta_0\right) \\ \log\Lambda_n\left(\theta_n, \theta_0\right) \end{pmatrix} \to N_{k+1}\left(\begin{pmatrix} -b(\delta) \\ -\frac{1}{2}\delta^{\mathrm{T}} J(\theta_0)\delta \end{pmatrix}, \begin{pmatrix} D\Sigma D^{\mathrm{T}} & -b(\delta) \\ -b(\delta)^{\mathrm{T}} & \delta^{\mathrm{T}} J(\theta_0)\delta \end{pmatrix}\right) \text{ in distribution.}$$

Invoke Cramer-Wold device again to arrive at the desired statement (4). The proof is thus completed. □

COROLLARY 1. *Given $\mathscr{F}_{t-1}$, suppose $\hat{\theta}$ denotes conditional maximum likelihood estimator of $\theta$ based on the conditional likelihood $\prod_{t=1}^{n} f_\theta(x_t|\mathscr{F}_{t-1})$. Then under A1-A4 and A6, and under local misspecification,*

$$n^{1/2}(\hat{\theta} - \theta_0) \to N_k\left(\delta, J(\theta_0)^{-1}\right) \text{ in distribution}$$

*Proof.* Given $\mathscr{F}_{t-1}$, conditional maximum likelihood estimator is an $M$-estimator with $\psi(X_t, \theta) = S(X_t, \theta)$. Moreover A6 implies that A5 with $J(\theta_0) = -d\phi(\theta_0) = \Sigma = E$, cf. van der Vaart (2000), Lemma 5.4. The proof is now immediate due to Theorem 2. □

*Remark* 1. Corollary 1 is seen to be a generalization of Theorem 5.1 of Claeskens and Hjort (2008) to stationary and strong mixing stochastic processes. Moreover, we obtain it as a special case of asymptotic normality of $M$-estimators with weaker regularity assumptions on likelihood than those of Claeskens and Hjort.

## 4. THE FOCUSED INFORMATION CRITERION

### 4·1. *Partitioned matrices*

Without loss of generality, suppose the parameter space $\Theta \subseteq \mathbb{R}^k$ is partitioned in two subspaces $A \subseteq \mathbb{R}^p$ and $B \subseteq \mathbb{R}^q$, $p + q = k$. Accordingly, let the labeling parameter vector $\theta$ be partitioned as $\theta = (\alpha, \beta), \alpha \in A, \beta \in B$. Let $P_{\alpha,\beta}$ denote the probability measure corresponding to the widest, full model for the data. We assume that there exists a unique $\beta_0$ in the interior of $B$ such that setting $\beta = \beta_0$ gives rise to the smallest, narrow model for the data. For example, the Weibull or gamma density with unknown scale parameter $\alpha$ and shape parameter $\beta$ is a wide model while if $\beta$ is set equal to 1, this model reduces to the exponential model which is a narrow model in this case. All the models that lie between the wide and narrow model are called submodels. We assume that all submodels retain $\alpha$-part of the wide model while elements in $\beta$ vector are set equal to the corresponding elements in $\beta_0$ vector.

Set $\psi_t = \psi(X_t, \theta)$ which we partition according to $\alpha$ and $\beta$ components as $\psi_t = (\psi_{t,\alpha}, \psi_{t,\beta})^{\mathrm{T}}$. Such partitioning is possible since it is assumed implicitly that function $\psi$ appears as a gradient of a criterion function to be optimized. Thus, the two components of vector $\psi_t$ denote the gradients of the criterion function with respect to $\alpha$ and $\beta$ respectively. We index the submodels with a set $\mathcal{S} \subseteq \{1, 2, ..., q\}$. Let $|\mathcal{S}|$ denote the cardinality of $\mathcal{S}$. Let $\pi_{\mathcal{S}}$ denote $|\mathcal{S}| \times q$ denote projection matrix such that for any $v \in \mathbb{R}^q$, $v_{\mathcal{S}} = \pi_{\mathcal{S}} v$ denotes the projection of $v$ from $\mathbb{R}^q$ into $\mathbb{R}^{\mathcal{S}}$. Since all the submodels retain the $\alpha$-part, the local misspecification involves only the perturbations around $\beta$. Thus, $\delta$ can be partitioned into $(0_p, \eta)^{\mathrm{T}}$ where $\eta \in \mathbb{R}^q$. Consequently, $P_{n,\alpha,\beta_0+n^{-1/2}\eta}$ denotes a sequence of locally misspecified probability measures.

Corresponding to a submodel $\mathcal{S}$, let $d_{\mathcal{S}} = d_{\mathcal{S}}(\phi)$, $D_{\mathcal{S}}$, $\Sigma_{\mathcal{S}}$, $E_{\mathcal{S}}$ denote the quantities $d\phi, D, \Sigma, E$ in Sections 2. It is understood that in each of these partitions, the components of $\beta$ that do not appear in $\mathcal{S}$ are set equal to their corresponding null value in $\beta_0$.

It can be seen that

$$\psi_{t,\mathcal{S}} = \begin{pmatrix} \psi_{t,\alpha} \\ \pi_{\mathcal{S}} \psi_{t,\beta} \end{pmatrix}.$$

Thus,

$$d_{\mathcal{S}} = \begin{pmatrix} d_{\alpha,\alpha} & d_{\alpha,\beta}\pi_{\mathcal{S}}^{\mathrm{T}} \\ \pi_{\mathcal{S}} d_{\beta,\alpha} & \pi_{\mathcal{S}} d_{\beta,\beta}\pi_{\mathcal{S}}^{\mathrm{T}} \end{pmatrix}$$

and

$$\Sigma_{\mathcal{S}} = \begin{pmatrix} \Sigma_{\alpha,\alpha} & \Sigma_{\alpha,\beta}\pi_{\mathcal{S}}^{\mathrm{T}} \\ \pi_{\mathcal{S}}\Sigma_{\beta,\alpha} & \pi_{\mathcal{S}}\Sigma_{\beta,\beta}\pi_{\mathcal{S}}^{\mathrm{T}} \end{pmatrix}.$$

Let $D_{\mathcal{S}}$ denote $d_{\mathcal{S}}^{-1}$ which can be partitioned accordingly as

$$D_{\mathcal{S}} = \begin{pmatrix} d_{\mathcal{S}}^{\alpha,\alpha} & d_{\mathcal{S}}^{\alpha,\beta} \\ d_{\mathcal{S}}^{\beta,\alpha} & d_{\mathcal{S}}^{\beta,\beta} \end{pmatrix}.$$

Since all the submodels retain $\alpha$ component, the submatrices $d_{\mathcal{S}}^{\alpha,\alpha}, d_{\mathcal{S}}^{\alpha,\beta}$ and $d_{\mathcal{S}}^{\beta,\beta}$ are of dimensions $p \times p, p \times |\mathcal{S}|$ and $|\mathcal{S}| \times |\mathcal{S}|$ respectively.

The score vector $S_t = S\left(X_t, \theta\right)$ can be partitioned similarly.

$$S_{t,\mathcal{S}} = \begin{pmatrix} S_{t,\alpha} \\ \pi_{\mathcal{S}} S_{t,\beta} \end{pmatrix}.$$

Finally, partition the matrix $E_{\mathcal{S}}$ as

$$E_{\mathcal{S}} = \begin{pmatrix} E_{\alpha,\alpha} & E_{\alpha,\beta}\pi_{\mathcal{S}}^{\mathrm{T}} \\ \pi_{\mathcal{S}}E_{\beta,\alpha} & \pi_{\mathcal{S}}E_{\beta,\beta}\pi_{\mathcal{S}}^{\mathrm{T}} \end{pmatrix},$$

where, $E_{\alpha,\alpha} = E\left(\psi_{0,\alpha} S_{0,\alpha}^{\mathrm{T}}\right)$ etc. Suppose $E$ is non-singular. Then, $E^{-1}$ is partitioned as

$$E^{-1} = \begin{pmatrix} E^{\alpha,\alpha} & E^{\alpha,\beta} \\ E^{\beta,\alpha} & E^{\beta,\beta} \end{pmatrix}.$$

### 4·2. *Asymptotic normality under subset models*

Let $\hat{\alpha}_{\mathcal{S}}$ and $\hat{\beta}_{\mathcal{S}}$ denote respectively estimators of $\alpha$ and $\beta_{\mathcal{S}}$ under model $\mathcal{S}$.

THEOREM 3. *Under assumptions of Theorem 2,*

$$n^{1/2}\begin{pmatrix} \hat{\alpha}_{\mathcal{S}} - \alpha_0 \\ \hat{\beta}_{\mathcal{S}} - \beta_{0,\mathcal{S}} \end{pmatrix} = \begin{pmatrix} d_{\mathcal{S}}^{\alpha,\alpha}E_{\alpha,\beta}\eta + d_{\mathcal{S}}^{\alpha,\alpha}N_\alpha + d_{\mathcal{S}}^{\alpha,\beta}\pi_{\mathcal{S}}E_{\beta,\beta}\eta + d_{\mathcal{S}}^{\alpha,\beta}N_{\beta_{\mathcal{S}}} \\ d_{\mathcal{S}}^{\beta,\alpha}E_{\alpha,\beta}\eta + d_{\mathcal{S}}^{\beta,\alpha}N_\alpha + d_{\mathcal{S}}^{\beta,\beta}\pi_{\mathcal{S}}E_{\beta,\beta}\eta + d_{\mathcal{S}}^{\beta,\beta}N_{\beta_{\mathcal{S}}} \end{pmatrix} + o_{P_{\alpha_0,\beta_{0,\mathcal{S}}}}(1),$$

*where* $N_\alpha \sim N_p\left(0_p, \Sigma_{\alpha,\alpha}\right)$ *and* $N_{\beta_{\mathcal{S}}} \sim N_{|\mathcal{S}|}\left(0_{|\mathcal{S}|}, \pi_{\mathcal{S}}\Sigma_{\beta,\beta}\pi_{\mathcal{S}}^{\mathrm{T}}\right)$.

*Proof.* Identify $T$ in Theorem A1 in Appendix with $M$-functional (1). Then, Lemma 1 along with equations (A4)-(A6) in Appendix implies that $n^{1/2}\begin{pmatrix} \hat{\alpha} - \alpha_0 \\ \hat{\beta} - \beta_0 \end{pmatrix}$ is asymptotically equivalent to $-n^{-1/2}D\sum_{t=1}^{n}\psi_t$ under $P_{\alpha_0,\beta_0}$. Applying Slutsky's theorem to Theorem A1 hence implies

$$n^{-1/2}\sum_{t=1}^{n}\psi_t \rightarrow N_{p+q}\left(0, \Sigma\right) \text{ in distribution.} \tag{8}$$

Le Cam's third lemma in view of equations (6) and (8) yields that under the local misspecification,

$$n^{-1/2}\sum_{t=1}^{n}\psi_t \rightarrow N_k\left(E\delta, \Sigma\right) \text{ in distribution.} \tag{9}$$

But,

$$E\delta = \begin{pmatrix} E_{\alpha,\beta}\eta \\ E_{\beta,\beta}\eta \end{pmatrix}.$$

Applying Cramer-Wold device we have,

$$n^{-1/2}\sum_{t=1}^{n}\psi_{t,\mathcal{S}}\rightarrow\begin{pmatrix}E_{\alpha,\beta}\eta+N_{\alpha}\\ \pi_{\mathcal{S}}E_{\beta,\beta}\eta+N_{\beta_{\mathcal{S}}}\end{pmatrix}\text{ in distribution.}$$

Since $n^{1/2}\begin{pmatrix}\hat{\alpha}_{\mathcal{S}}-\alpha_0\\ \hat{\beta}_{\mathcal{S}}-\beta_{0,\mathcal{S}}\end{pmatrix}$ is asymptotically equivalent to $n^{-1/2}\sum_{t=1}^{n}\psi_{t,\mathcal{S}}$, the proof now follows by virtue of Slutsky's lemma. □

### 4·3. *Derivation of* FIC

Fix $t\geq 2$. Given $\mathscr{F}_{t-1}$, let $\mu_{t-1}=\mu\left(\alpha,\beta,\mathscr{F}_{t-1}\right)$ denote the real valued differentiable focus function of parameters. For instance, the focus function while selecting the order of an autoregressive process can be the conditional mean of the process given the lag of order $r$. For any submodel $\mathcal{S}$, let the true value of the focus under the local misspecification be given as $\mu_{t-1}^{\mathcal{S}}=\mu\left(\alpha_0,\beta_{0,\mathcal{S}}+n^{-1/2}\pi_{\mathcal{S}}\eta,\mathscr{F}_{t-1}\right)$. Let $\hat{\mu}_{t-1}^{\mathcal{S}}=\mu\left(\hat{\alpha_{\mathcal{S}}},\hat{\beta_{\mathcal{S}}},\mathscr{F}_{t-1}\right)$ denote an estimator of $\mu_{t-1}$ under a submodel $\mathcal{S}$. In order to derive the FIC for submodel $\mathcal{S}$, asymptotic normal distribution of $n^{1/2}\left(\hat{\mu}_{t-1}^{\mathcal{S}}-\mu_{t-1}^{\mathcal{S}}\right)$ needs to be specified completely.

Let $\zeta_{\mathcal{S}}=(\zeta_{\mathcal{S},\alpha},\zeta_{\mathcal{S},\beta})$ and $\mathcal{C}_{\mathcal{S}}=\begin{pmatrix}\mathcal{C}_{\mathcal{S},\alpha,\alpha}&\mathcal{C}_{\mathcal{S},\alpha,\beta}\\ \mathcal{C}_{\mathcal{S},\beta,\alpha}&\mathcal{C}_{\mathcal{S},\beta,\beta}\end{pmatrix}$ denote respectively the mean vector and covariance matrix of asymptotic normal distribution of $n^{1/2}\begin{pmatrix}\hat{\alpha}_{\mathcal{S}}-\alpha_0\\ \hat{\beta}_{\mathcal{S}}-\beta_{0,\mathcal{S}}\end{pmatrix}$. By means of the Theorem 3, it can be deduced that

$$\begin{aligned}\zeta_{\mathcal{S},\alpha}&=W_{\mathcal{S}}^{1}\eta,\\ \zeta_{\mathcal{S},\beta}&=W_{\mathcal{S}}^{2}\eta,\end{aligned}$$

where $W_{\mathcal{S}}^{1}$ and $W_{\mathcal{S}}^{2}$ are matrices of dimensions $p\times q$ and $|\mathcal{S}|\times q$ respectively given as

$$\begin{aligned}W_{\mathcal{S}}^{1}&=d_{\mathcal{S}}^{\alpha,\alpha}E_{\alpha,\beta}+d_{\mathcal{S}}^{\alpha,\beta}\pi_{\mathcal{S}}E_{\beta,\beta},\\ W_{\mathcal{S}}^{2}&=d_{\mathcal{S}}^{\beta,\alpha}E_{\alpha,\beta}+d_{\mathcal{S}}^{\beta,\beta}\pi_{\mathcal{S}}E_{\beta,\beta}\end{aligned}.$$

Moreover,

$$\begin{aligned}\mathcal{C}_{\mathcal{S},\alpha,\alpha}&=d_{\mathcal{S}}^{\alpha,\beta}\pi_{\mathcal{S}}\Sigma_{\beta,\beta}\pi_{\mathcal{S}}^{\mathrm{T}}d_{\mathcal{S}}^{\beta,\alpha}+2d_{\mathcal{S}}^{\alpha,\beta}\pi_{\mathcal{S}}\Sigma_{\beta,\alpha}d_{\mathcal{S}}^{\alpha,\alpha}+d_{\mathcal{S}}^{\alpha,\alpha}\Sigma_{\alpha,\alpha}d_{\mathcal{S}}^{\alpha,\alpha}\\ \mathcal{C}_{\mathcal{S},\alpha,\beta}&=d_{\mathcal{S}}^{\alpha,\alpha}\Sigma_{\alpha,\alpha}d_{\mathcal{S}}^{\alpha,\beta}+d_{\mathcal{S}}^{\alpha,\beta}\pi_{\mathcal{S}}\Sigma_{\beta,\beta}\pi_{\mathcal{S}}^{\mathrm{T}}d_{\mathcal{S}}^{\beta,\beta}+d_{\mathcal{S}}^{\alpha,\alpha}\Sigma_{\alpha,\beta}\pi_{\mathcal{S}}^{\mathrm{T}}d_{\mathcal{S}}^{\beta,\beta}+d_{\mathcal{S}}^{\alpha,\beta}\pi_{\mathcal{S}}\Sigma_{\beta,\alpha}d_{\mathcal{S}}^{\beta,\alpha}\\ \mathcal{C}_{\mathcal{S},\beta,\beta}&=d_{\mathcal{S}}^{\beta,\beta}\pi_{\mathcal{S}}\Sigma_{\beta,\beta}\pi_{\mathcal{S}}^{\mathrm{T}}d_{\mathcal{S}}^{\beta,\beta}+2d_{\mathcal{S}}^{\beta,\alpha}\Sigma_{\alpha,\beta}\pi_{\mathcal{S}}^{\mathrm{T}}d_{\mathcal{S}}^{\beta,\beta}+d_{\mathcal{S}}^{\beta,\alpha}\Sigma_{\alpha,\alpha}d_{\mathcal{S}}^{\alpha,\beta}.\end{aligned}$$

Suppose $\mathcal{I}_q$ denotes an identity matrix of order $q$.
Define

$$m_{\mathcal{S}}=\left(\frac{\partial\mu_{t-1}}{\partial\alpha_0}\right)^{\mathrm{T}}W_{\mathcal{S}}^{1}+\left(\frac{\partial\mu_{t-1}}{\partial\beta_0}\right)^{\mathrm{T}}\left(\pi_{\mathcal{S}}^{\mathrm{T}}W_{\mathcal{S}}^{2}-\mathcal{I}_q\right)\tag{10}$$

and

$$\begin{aligned}V_{\mathcal{S}}=&\left(\frac{\partial\mu_{t-1}}{\partial\alpha_0}\right)^{\mathrm{T}}\mathcal{C}_{\mathcal{S},\alpha,\alpha}\left(\frac{\partial\mu_{t-1}}{\partial\alpha_0}\right)+2\left(\frac{\partial\mu_{t-1}}{\partial\alpha_0}\right)^{\mathrm{T}}\mathcal{C}_{\mathcal{S},\alpha,\beta}\pi_{\mathcal{S}}\left(\frac{\partial\mu_{t-1}}{\partial\beta_0}\right)\\ &+\left(\frac{\partial\mu_{t-1}}{\partial\beta_0}\right)^{\mathrm{T}}\pi_{\mathcal{S}}^{\mathrm{T}}\mathcal{C}_{\mathcal{S},\beta,\beta}\pi_{\mathcal{S}}\left(\frac{\partial\mu_{t-1}}{\partial\beta_0}\right).\end{aligned}\tag{11}$$

The following corollary to Theorem 3 is instrumental in deriving FIC expressions.

COROLLARY 2. *The asymptotic mean squared error of* $n^{1/2}\left(\hat{\mu}_{t-1}^{\mathcal{S}}-\mu_{t-1}^{\mathcal{S}}\right)$ *is* $m_{\mathcal{S}}^{\mathrm{T}}\eta\eta^{\mathrm{T}}m_{\mathcal{S}}+V_{\mathcal{S}}$*, where* $m_{\mathcal{S}}$ *and* $V_{\mathcal{S}}$ *are as in (10) and (11) respectively.*

*Proof.* Since $\mu_{t-1}$ is differentiable focus function, first order Taylor's expansion around $\beta_0$ yields

$$n^{1/2}\left(\hat{\mu}_{t-1}^{\mathcal{S}}-\mu_{t-1}^{\mathcal{S}}\right)=n^{1/2}\left(\hat{\mu}_{t-1}^{\mathcal{S}}-\mu_{t-1}\left(\alpha_0,\beta_{0,\mathcal{S}}\right)\right)-\eta^{\mathrm{T}}\pi_{\mathcal{S}}^{\mathrm{T}}\pi_{\mathcal{S}}\left(\frac{\partial\mu_{t-1}}{\partial\beta_0}\right).$$

Delta method applied to the right hand side of the above relation coupled with Theorem 3 shows that

$$n^{1/2}\left(\hat{\mu}_{t-1}^{\mathcal{S}}-\mu_{t-1}^{\mathcal{S}}\right)\rightarrow N_1\left(m_{\mathcal{S}}^{\mathrm{T}}\eta,V_{\mathcal{S}}\right)\text{ in distribution.}$$

The proof now follows immediately. □

An estimate of this mean squared error based on the data at one's disposal is the FIC score for submodel $\mathcal{S}$. The quantities $m_{\mathcal{S}}$ and $V_{\mathcal{S}}$ can be estimated consistently by plugging the estimates of parameters under wide model in (10) and (11) respectively. It is also admissible to estimate these quantities under the narrow model, in which case, $\beta$ is replaced by $\beta_0$ and $\alpha$ is replaced by its estimate under narrow model.

*Remark* 2. In the case of the wide (largest) model, the projection matrix is an identity matrix of order $q$ while for the narrow (smallest) model, it is a null matrix of order $q$. Since the narrow model does not involve an unknown $\beta$ component, the matrices in (10) and (11) are obtained using Theorem 3 following the proof of Corollary 2. Hence,

$$m_{\text{narrow}}=\left(\frac{\partial\mu_{t-1}}{\partial\alpha_0}\right)^{\mathrm{T}}d_{\alpha,\alpha}^{-1}E_{\alpha,\beta},$$

$$V_{\text{narrow}}=\left(\frac{\partial\mu_{t-1}}{\partial\alpha_0}\right)^{\mathrm{T}}d_{\alpha,\alpha}^{-1}\Sigma_{\alpha,\alpha}d_{\alpha,\alpha}^{-1}\left(\frac{\partial\mu_{t-1}}{\partial\alpha_0}\right).$$

Identifying each of the matrices $d$, $E$ and $\Sigma$ to be the Fisher information matrix, the Corollary 2 above implies Corollary 5.1 of Claeskens and Hjort (2008) when the focus is estimated using maximum likelihood estimators.

### 4·4. *Estimation of* $\eta\eta^{\mathrm{T}}$

Estimation of local misspecification radius $\eta$ turns out to be elusive due to non-existence of any consistent estimator for the same. However, it is an estimator of $\eta\eta^{\mathrm{T}}$ that is indeed needed while estimating the squared bias component of FIC. It is possible to device plug-in as well as asymptotically unbiased estimators of $\eta\eta^{\mathrm{T}}$ that work for practical purposes.

Cramer-Wold device applied to (9) implies

$$G_n\rightarrow N_k\left(\eta,H\right)\text{ in distribution,}\tag{12}$$

where

$$G_n=n^{-1/2}\left(E^{\beta,\alpha}\sum_{t=1}^{n}\psi_{t,\alpha}+E^{\beta,\beta}\sum_{t=1}^{n}\psi_{t,\beta}\right)$$

and

$$H=E^{\beta,\alpha}\left(\Sigma_{\alpha,\alpha}E^{\alpha,\beta}+\Sigma_{\alpha,\beta}E^{\beta,\beta}\right)+E^{\beta,\beta}\left(\Sigma_{\beta,\alpha}E^{\alpha,\beta}+\Sigma_{\beta,\beta}E^{\beta,\beta}\right).$$

Equation (12) implies that $G_n$ is asymptotically unbiased for $\eta$. Hence, a plug-in estimator of $\eta\eta^T$ is a sample version of $G_n G_n^{\mathrm{T}}$.

Relation (12) also implies

$$E\left(G_n G_n^{\mathrm{T}}\right) \to \eta\eta^{\mathrm{T}} + H.$$

Hence, the bias-corrected version of the above plug-in estimator is a sample version $G_n G_n^{\mathrm{T}} - H$.

*Remark* 3. It is clear that for any $\eta \in \mathbb{R}^q$, $\eta\eta^{\mathrm{T}}$ is a $q \times q$ positive semidefinite matrix. The bias-corrected version, however, may not be positive semidefinite in general and thus it may not be a valid estimator. On the other hand, although a plug-in estimator is biased, it is positive semidefinite with probability 1. Hence, albeit biased, a plug-in estimator may be preferred over an asymptotically unbiased estimator. Lohmeyer et al. (2018) give an account of the impacts of both these estimators on performance of FIC.

*Remark* 4. The above formulae of FIC are derived for a fixed time point $t$. In practice, the data are observed on several time-points. Thus, the focus function indeed may be time-varying in general. The final FIC score for model selection is then computed by any suitable aggregation of FIC scores at individual time-points such as the average of individual scores as suggested by Rohan and Ramanathan (2011).

## 5. Important special cases

### 5·1. *Efficient estimators based on conditional likelihood*

Given $\mathcal{F}_{t-1}$, let the estimating function $\psi_t$ denote the conditional score vector $S_t$ . Then, $E = \Sigma = J(\theta_0)$ and $D_{\mathcal{S}} = J_{\mathcal{S}}$, the Fisher information matrix for submodel $\mathcal{S}$. Partition $J(\theta_0)$ as well as $J_{\mathcal{S}}$ according to $\alpha$ and $\beta$ components as demonstrated in Subsection 4.1. Set

$$\omega = J_{\beta,\alpha} J_{\alpha,\alpha}^{-1} \frac{\partial \mu_{t-1}}{\partial \alpha} - \frac{\partial \mu_{t-1}}{\partial \beta}.$$

$$\tau^2 = \left(\frac{\partial \mu_{t-1}}{\partial \alpha}\right)^{\mathrm{T}} J_{\alpha,\alpha}^{-1} \frac{\partial \mu_{t-1}}{\partial \alpha}.$$

Set $\omega_{\mathcal{S}} = \pi_{\mathcal{S}}\omega, Q = J^{\beta,\beta}$, $Q_{\mathcal{S}} = J_{\mathcal{S}}^{\beta,\beta}$ and $\Gamma_{\mathcal{S}} = \pi_{\mathcal{S}}^{\mathrm{T}} Q_{\mathcal{S}} \pi_{\mathcal{S}} Q^{-1}$. By means of the formulae for inverse of partitioned matrices in Harville (1997), Subsection 8.5, a few algebraic efforts lead to the simplified expressions $m_{\mathcal{S}}$ and $V_{\mathcal{S}}$ as below.

$$m_{\mathcal{S}} = \omega^{\mathrm{T}}\left(I_q - \Gamma_{\mathcal{S}}\right).$$
$$V_{\mathcal{S}} = \tau^2 + \omega^{\mathrm{T}} Q_{\mathcal{S}}\omega.$$

With these simplifications, our formula for asymptotic mean squared error corresponds to that of Claeskens and Hjort (2008), Equation 6.6. Since efficient estimators are asymptotically linear in score vector, asymptotically unbiased estimator of $\eta\eta^{\mathrm{T}}$ is given $D_n D_n^{\mathrm{T}} - Q$ where $D_n = n^{1/2}\left(\hat{\beta} - \beta_0\right)$. Thus, our general expression of FIC reduces to the simpler form as given by Claeskens and Hjort (2008), Equation 6.7. Following is a further simpler form of FIC provided by Claeskens and Hjort (2008), Equation 6.1.

$$\text{FIC}_\mu(\mathcal{S}) = \left(\hat{\omega}^{\mathrm{T}} D_n - \hat{\omega}^{\mathrm{T}} \Gamma_{\mathcal{S}} D_n\right)^2 + 2\hat{\omega}_{\mathcal{S}}^{\mathrm{T}} \hat{Q}_{\mathcal{S}} \hat{\omega}_{\mathcal{S}}. \tag{13}$$

A noteworthy feature of expression (13) is that it does not require computation of $\tau^2$. Moreover, the above formula is applicable for model selection in the dependent setup as well since it is now viewed as an immediate consequence and a special case of our general theory.

### 5·2. *Conditional least square estimators*

Let $k$ denote a nonnegative integer. Set $\mathscr{F}_{t-1} = \sigma\left(X_{t-1}, .., X_{t-k}\right), t \geq k$. Define a real valued measurable function $g$ as

$$g\left(\theta, \mathscr{F}_{t-1}\right) = E_\theta\left(X_t | \mathscr{F}_{t-1}\right).$$

The conditional least square estimator of $\theta$ as proposed by Klimko and Nelson (1978) is obtained by minimizing the residual sum of squares

$$\sum_{t=k+1}^{n} \left(x_t - g\left(\theta, \mathscr{F}_{t-1}\right)\right)^2.$$

Suppose $g$ is a differentiable function of $\theta$. Let $\nabla g\left(\theta, \mathscr{F}_{t-1}\right)$ denote the gradient of $g$ with respect to $\theta$. Assume that $\nabla g\left(\theta_0, \mathscr{F}_{t-1}\right)$ is non-vanishing at some $\theta_0$ in the interior of $\Theta$. Then for any $t \geq k$,

$$\psi_t = \nabla g\left(\theta, \mathscr{F}_{t-1}\right)\left(x_t - g\left(\theta, \mathscr{F}_{t-1}\right)\right).$$

Let $\{\epsilon_t, t \geq k\}$ denote a martingale difference sequence where $\epsilon_t = X_t - g\left(\theta, \mathscr{F}_{t-1}\right)$. In view of Theorems 2.2 and 3.2 of Klimko and Nelson, the matrices $D$ and $\Sigma$ in our Theorem 2 are identified as

$$D = \left(E\left(\nabla g\left(\theta_0, \mathscr{F}_{k-1}\right)\left(\nabla g\left(\theta_0, \mathscr{F}_{k-1}\right)\right)^{\mathrm{T}}\right)\right)^{-1},$$

and

$$\Sigma = E\left(\epsilon_k^2 \nabla g\left(\theta_0, \mathscr{F}_{k-1}\right)\left(\nabla g\left(\theta_0, \mathscr{F}_{k-1}\right)\right)^{\mathrm{T}}\right).$$

Let $f(\epsilon_t)$ denote the error density. Let

$$S_t(\theta) = \nabla \log f\left(x_t - g\left(\theta, \mathscr{F}_{k-1}\right)\right)$$

denote the score vector. The matrix $E$ is given as

$$E = E\left(\epsilon_k \nabla g\left(\theta_0, \mathscr{F}_{k-1}\right) S_k(\theta_0)^{\mathrm{T}}\right).$$

Substituting these expressions in (10) and (11) yields the desired FIC expressions. The sample versions of the above expectations consistently estimate the matrices $D, \Sigma$ and $E$ owing to the stationarity and ergodicity of the underlying process.

*Remark* 5. The common model selection problems in the time series models involve selection of suitable $k$ reflecting the order of time-dependence. The focus function in FIC can thus be chosen as the conditional expectation $g$. Linear parametric function $g$ corresponds to the well-known autoregressive time series models. The FIC for Gaussian linear time series models is proposed by Rohan and Ramanathan (2011). It is now clear that this forms the special case of our general setup.

### 5·3. *Method of moments*

For the sake of simplicity, assume that $\{X_n\}$ denotes a sequence of independent and identically distributed random variables on $\left(\Omega, \mathscr{F}, P_\theta\right), \theta \in \mathbb{R}^k$. Choose a measurable function $h : (\Omega, \mathscr{F}) \to \left(\mathbb{R}^k, \mathscr{B}(\mathbb{R}^k)\right)$, where $\mathscr{B}(\mathbb{R}^k)$ is a $k$-dimensional Borel $\sigma$-algebra, such that $E||h(X_1)|| < \infty$. Denote by $\nu(\theta)$ a $k \times 1$ mean vector $E\left(h(X_1)\right)$. The method of moments estimator is an $M$-estimator with estimating function in (2) as $\psi_i = h(X_i) - \nu(\theta), i = 1, 2, ..., n$. Suppose that $\nu(\theta)$ is differentiable with $\nabla\nu(\theta_0)$ as the locally homeomorphic derivative at $\theta_0$. It is easy to see that $D = \nabla\nu(\theta_0)^{-1}$, and $\Sigma = E\left(h(X_1)h(X_1)^{\mathrm{T}}\right)$. Let $S_1 = \frac{\partial \log f(x_1,\theta)}{\partial\theta_0}$ denote the score vector evaluated at $\theta_0$. Since $E(S_1) = 0_k$, the matrix $E$ is same as $E\left(h(X_1)S_1^{\mathrm{T}}\right)$. Partition each of these matrices according to the components of the narrow and the wide model as in Section 4. The computation of FIC then follows from Corollary 2.

## 6. Application: Asymmetric Laplace Autoregressive Models

### 6·1. *Model and problem of model selection*

The asymmetric version of Laplace (Double exponential) models is proposed by Kotz et al. (2001) with applications to financial modeling. The autoregressive models of order $q \geq 1$ driven by asymmetric Laplace innovations are studied in details by Trindale et al. (2009).

A stochastic process $\{X_t, t \in \mathbb{Z}^+\}$ is said to be an autoregressive process of order $q$ driven by asymmetric Laplace noise if it is a stationary solution of the equations

$$X_t = \alpha_1 X_{t-1} + ... + \alpha_p X_{t-p} + \alpha_{p+1} X_{t-p-1} + ... + \alpha_q X_{t-q} + Z_t, \tag{14}$$

where $\{Z_t, t \geq q\}$ is a sequence of independent and identically distributed asymmetric Laplace random variables with location $\theta \in \mathbb{R}$, scale $\sigma > 0$ and skewness parameter $\kappa \in \mathbb{R}$. The probability density of asymmetric Laplace random variables is given by Kotz et al. (2001) as

$$\begin{aligned} f(y, \theta, \sigma, \kappa) &= \frac{2^{1/2}\kappa}{\sigma(1+\kappa^2)} exp\{\frac{2^{1/2}}{\kappa\sigma}(y-\theta)\} \text{ if } y \leq \theta \\ &= \frac{2^{1/2}\kappa}{\sigma(1+\kappa^2)} exp\{-\frac{2^{1/2}\kappa}{\sigma}(y-\theta)\} \text{ if } y > \theta. \end{aligned}$$

As suggested by Trindale et al. (2009), we set

$$\theta = 2^{-1/2}\sigma(\kappa - 1/\kappa),$$

so that $E(Z_t) = 0$ for every $t$ which is a typical assumption of autoregressive models. Henceforth, we denote such a process by $ALAR(q)$.

The selection problem in models given by (14) is two-fold. One aspect of model selection here is determination of the order of dependence of the process. The other aspect is the selection between skewed and classical, symmetric Laplace density to model the innovations. Both these problems need to be tackled simultaneously since both these aspects are likely to be intertwined and interdependent. The focused model selection here involves the simultaneous selection of order and symmetry of error density in the local misspecification identifying a suitable parametric focus function. The three natural choices of foci are conditional mean of the process, coefficient of skewness of error density and a linear combination of these two foci. We derive the FIC to tackle these selection problems using conditional maximum likelihood estimator. The derivation for FIC based on conditional least square estimators follows on the same lines and is omitted to save the space.

The model given by (14) can be regarded as the richest model with $q+2$ parameters. Set $\alpha = (\alpha_1, ..., \alpha_p)$, $\gamma = (\alpha_{p+1}, ..., \alpha_q) = (\gamma_1, ..., \gamma_{q-p})$ with $p < q$. Let $\lambda = (\alpha, \gamma, \kappa, \sigma)$ denote the full parameter vector of length $q+2$. Assume that the parameters $\gamma$ and $\kappa$ are subject to the perturbations in local misspecification setting while $p+1$ parameters $(\alpha, \sigma)$ are common in all the submodels. Setting $\kappa = 1$ leads to autoregressive processes of order $q$ driven by Laplace error abbreviated as $LAR(q)$ process. For any $0 \leq p < q$, the sequence $LAR(p) \subseteq ALAR(p) \subset LAR(q) \subset ALAR(q)$ forms a nested sequence of models with dimensions $p+1, p+2, q+1$ and $q+2$ respectively. The wide model $ALAR(q)$ reduces to the narrow model $LAR(p)$ when $\gamma = 0_{q-p}$ and $\kappa = 1$. Hence, we assume that the true error density is locally misspecified as $f_{\text{true}}(z_t) = f(z_t, \alpha_0, n^{-1/2}\delta_1, 1 + n^{-1/2}\delta_2, \sigma_0)$; $t \geq q$, $\delta_1 \in \mathbb{R}^{q-p}$ and $\delta_2 \in \mathbb{R}$.

### 6·2. *The* FIC

Suppose that the roots of the characteristic polynomial of degree $q$ associated with the autoregressive process (14) lie outside the unit circle, so that the process is stationary and ergodic. Suppose $\hat{\lambda}$ denotes the conditional maximum likelihood estimator of $\lambda$ under the wide model (14). Let $J_{wide}(\lambda)$ denote the Fisher information matrix. Since asymmetric Laplace density is a member of location-scale family and is differentiable with respect to location $\theta$ except on the set of Lebesgue measure 0, it is $\mathcal{L}^2$-differentiable (Le Cam and Yang, 2012, Lemma 1). Thus, in view of our Corollary 1, it can be claimed that

$$n^{1/2}\begin{pmatrix} \hat{\alpha} - \alpha_0 \\ \hat{\gamma} \\ \hat{\kappa} - 1 \\ \hat{\sigma} - \sigma_0 \end{pmatrix} \to N_{q+2}\left[\begin{pmatrix} 0_p \\ \delta_1 \\ \delta_2 \\ 0 \end{pmatrix}, J_{wide}^{-1}(\lambda)\right] \text{ in distribution.}$$

*Remark* 6. Although $\mathcal{L}^2$-differentiable, Laplace density is not Cramer-regular due to non-differentiability of the density at $\theta$. Thus, the theory of Claeskens and Hjort (2003, 2008) is inapplicable to derive FIC since it relies extensively on higher order differentiability assumptions on underlying density. This is the chief reason we choose to demonstrate our theory with this example. Our general asymptotic theory is presented with much weaker assumptions. The assumption of $\mathcal{L}^2$-differentiability in particular suffices to claim the local asymptotic normality as above and no additional assumptions on higher ordered derivatives are needed.

Following Trindale et al. (2009), the Fisher information matrix for wide model can be partitioned according to $\lambda = (\alpha, \gamma, \kappa, \sigma)$ and in that order as

$$\begin{aligned} J_{wide}(\lambda) &= \begin{pmatrix} \tilde{J}(\alpha,\gamma) & 0_{q\times 2} \\ 0_{2\times q} & J^*(\kappa,\sigma) \end{pmatrix} \\ &= \begin{pmatrix} \tilde{J}_{\alpha,\alpha} & \tilde{J}_{\alpha,\gamma} & 0_p & 0_p \\ \tilde{J}_{\gamma,\alpha} & \tilde{J}_{\gamma,\gamma} & 0_{(q-p)} & 0_{(q-p)} \\ 0_p^{\mathrm{T}} & 0_{(q-p)}^{\mathrm{T}} & J^*_{\kappa,\kappa} & J^*_{\kappa,\sigma} \\ 0_p^{\mathrm{T}} & 0_{(q-p)}^{\mathrm{T}} & J^*_{\sigma,\kappa} & J^*_{\sigma,\sigma} \end{pmatrix}. \end{aligned}$$

Here $\tilde{J}(\alpha,\gamma)$ is $q \times q$ positive definite matrix which consists of functions of asymmetric Laplace density and autocovariance functions of the underlying process (Li and McLeod 1988) and $J^*(\kappa,\sigma)$ is $2 \times 2$ Fisher information matrix corresponding to asymmetric Laplace density.

Let $\mathcal{S}_1$ and $\mathcal{S}_2$ denote respectively the submodels $LAR(q)$ and $ALAR(p)$. Let $\pi_{\mathcal{S}_1} = [I_{q-p} \mid 0_{q-p}]$ and $\pi_{\mathcal{S}_2} = [0_{q-p}^{\mathrm{T}} \mid 1]$ denote the projection matrices of orders $(q-p)\times(q-p+1)$ and $1 \times (q-p+1)$ respectively. Matrix $\omega$ in subsection 5.1 can thus be partitioned as $\omega(\mu) = (\omega_{\mathcal{S}_1}, \omega_{\mathcal{S}_2})^{\mathrm{T}}$ where

$$\omega_{\mathcal{S}_1} = \tilde{J}_{\gamma,\alpha}\tilde{J}_{\alpha,\alpha}^{-1}\frac{\partial\mu}{\partial\alpha} - \frac{\partial\mu}{\partial\gamma}$$

and

$$\omega_{\mathcal{S}_2} = J^*_{\kappa,\sigma} {J^*_{\sigma,\sigma}}^{-1} \frac{\partial \mu}{\partial \sigma} - \frac{\partial \mu}{\partial \kappa}.$$

Similarly,

$$Q = \begin{pmatrix} Q_{\mathcal{S}_1} & 0_{q-p} \\ 0^{\mathrm{T}}_{q-p} & Q_{\mathcal{S}_2} \end{pmatrix},$$

where, $Q_{\mathcal{S}_1} = \tilde{J}^{-1}_{\gamma,\gamma}$ and $Q_{\mathcal{S}_2} = {J^*_{\kappa,\kappa}}^{-1}$. Define $D_{n_1} = n^{1/2}\hat{\gamma}, D_{n_2} = n^{1/2}(\hat{\kappa} - 1)$. Using these quantities in formula (13) yields the following expressions of FIC for wide model $ALAR(q)$, subset models $\mathcal{S}_1$, $\mathcal{S}_2$ and narrow model $LAR(p)$.

$$\begin{aligned} \mathrm{FIC}_\mu(\text{wide}) &= 2(\hat{\omega}^{\mathrm{T}}_{\mathcal{S}_1} \hat{Q}_{\mathcal{S}_1} \hat{\omega}_{\mathcal{S}_1} + \hat{\omega}^2_{\mathcal{S}_2} \hat{Q}_{\mathcal{S}_2}) \\ \mathrm{FIC}_\mu(\mathcal{S}_1) &= (\hat{\omega}_{\mathcal{S}_2} D_{n_2})^2 + 2\hat{\omega}^{\mathrm{T}}_{\mathcal{S}_1} \hat{Q}_{\mathcal{S}_1} \hat{\omega}_{\mathcal{S}_1} \\ \mathrm{FIC}_\mu(\mathcal{S}_2) &= (\hat{\omega}^{\mathrm{T}}_{\mathcal{S}_1} D_{n_1})^2 + 2\hat{\omega}^2_{\mathcal{S}_2} \hat{Q}_{\mathcal{S}_2} \\ \mathrm{FIC}_\mu(\text{narrow}) &= (\hat{\omega}^{\mathrm{T}}_{\mathcal{S}_1} D_{n_1} + \hat{\omega}_{\mathcal{S}_2} D_{n_2})^2 \end{aligned} \tag{15}$$

### 6·3. *Focus on conditional mean*

A typical choice of focus function in order selection problems is $\mu_1 = E_\lambda(X_t | X_{t-q}, ..., X_{t-1})$. Set $Y_t^{(1)} = (X_{t-1}, ..., X_{t-p})$ and $Y_t^{(2)} = (X_{t-p-1}, ..., X_{t-q})$. Since errors are independent with mean 0, $\mu_1 = \alpha^{\mathrm{T}} Y_t^{(1)} + \gamma^{\mathrm{T}} Y_t^{(2)}$. Note that with this choice of focus, $\omega_{\mathcal{S}_2} = 0$ since focus does not depend on the parameters of asymmetric Laplace density. Thus, $\mathrm{FIC}_{\mu_1}(\mathcal{S}_1) = \mathrm{FIC}_{\mu_1}(\text{wide})$ and $\mathrm{FIC}_{\mu_1}(\mathcal{S}_2) = \mathrm{FIC}_{\mu_1}(\text{narrow})$. Hence, the problem of focused model selection is reduced to selection between $LAR(q)$ and $LAR(p)$ models or $ALAR(q)$ and $ALAR(p)$ . This is logical since conditional mean of the process takes into consideration only the time-dependence and may not provide any information about the properties of error density. Hence, equations (15) yield $\mathrm{FIC}_{\mu_1}(\mathcal{S}_1) = 2\hat{\omega}^{\mathrm{T}}_{\mathcal{S}_1} \hat{Q}_{\mathcal{S}_1} \hat{\omega}_{\mathcal{S}_1}$ and $\mathrm{FIC}_{\mu_1}(\text{narrow}) = (\hat{\omega}^{\mathrm{T}}_{\mathcal{S}_1} D_{n_1})^2$ where $\omega_{\mathcal{S}_1} = \tilde{J}_{\gamma,\alpha} \tilde{J}^{-1}_{\alpha,\alpha} Y_t^{(1)} - Y_t^{(2)}$.

### 6·4. *Focus on third central moment of error density*

The third central moment is a measure of skewness and hence may be used as a focus parameter when the selection between symmetric and asymmetric versions of the underlying probability distribution is intended. In this case, $\mu_2 = E(Z_q - E(Z_q))^3$. Since the focus does not depend on autoregression coefficients $(\alpha, \gamma)$ , $\omega_{\mathcal{S}_1} = 0_q$. Hence the problem is reduced to selection between $ALAR(p)$ and $LAR(p)$ models. This is justified since if the symmetry of the error distribution alone is the focus of model selection, the candidate models may be of the equal orders. In view of (15), $\mathrm{FIC}_{\mu_2}(\mathcal{S}_2) = 2\hat{\omega}^2_{\mathcal{S}_2} \hat{Q}_{\mathcal{S}_2} = \mathrm{FIC}_{\mu_2}(\text{wide})$ and $\mathrm{FIC}_{\mu_2}(\text{narrow}) = (\hat{\omega}_{\mathcal{S}_2} D_{n_2})^2 = \mathrm{FIC}_{\mu_2}(\mathcal{S}_1)$. It can be seen that $\mu_2 = 2^{-1/2}\sigma^3(1/\kappa^3 - \kappa^3)$. Hence, $\omega_{\mathcal{S}_2} = 2^{-1/2} 3\sigma^2 J^*_{\kappa,\sigma} {J^*_{\sigma,\sigma}}^{-1}(1/\kappa^3 - \kappa^3) + 2^{-1/2}\sigma^3(3/\kappa^4 + 3\kappa^2)$.

### 6·5. *Focus on linear combination*

An ideal choice of focus for simultaneous selection of order and symmetry is a linear combination of above two foci. For any two reals $l_1$ and $l_2$, let $\mu_3 = l_1\mu_1 + l_2\mu_2$ denote a focus function. $l_1$ and $l_2$ are chosen to reflect the relative importance of time-dependence and symmetry aspects of model selection. It can be seen that $\omega(\mu_3) = (l_1\omega_{\mathcal{S}_1}, l_2\omega_{\mathcal{S}_2})^{\mathrm{T}}$. Thus, for a subset model $\mathcal{S} \in \{\text{wide}, \mathcal{S}_1, \mathcal{S}_2\}$, $\mathrm{FIC}_{\mu_3}(\mathcal{S}) = l_1^2 \mathrm{FIC}_{\mu_1}(\mathcal{S}) + l_2^2 \mathrm{FIC}_{\mu_2}(\mathcal{S})$, while $\mathrm{FIC}_{\mu_3}(\text{narrow}) =$

$\left(l_1|\text{FIC}_{\mu_1}(\text{narrow})|^{1/2} + l_2|\text{FIC}_{\mu_2}(\text{narrow})|^{1/2}\right)^2$. It is evident that the focused model selection with focus on $\mu_3$ shares the advantages and features of the model selection based on $\mu_1$ and $\mu_2$.

## 7. SIMULATIONS

### 7·1. *Order and symmetry of errors in asymmetric Laplace autoregressive models*

In order to examine the performance of FIC as opposed to the traditional information criteria like AIC and BIC, we carryout a number of simulation exercises with two different choices of focus parameters viz., conditional mean of the process and the unconditional third central moment of model residuals. In both the setups, 100 samples each of sizes 500, 2000 and 3500 were simulated and the number of times each of the information criterion of interest selects the correct model was recorded. With the focus on conditional mean, each of the models from $ALAR(1)$ to $ALAR(6)$ was chosen to be the correct model to simulate from. The number of times out of 100 each of the three criteria select the correct model against the remaining competing models are recorded in Table 1 below. Overall impression is that the performance of FIC is comparable

Table 1. *Simulation results with focus on conditional mean*

| $n = 500$ | | | | | | | | | |
|---|---|---|---|---|---|---|---|---|---|
| True Order | Coefficients | | | | | | Scores | | |
| | $\alpha_1$ | $\alpha_2$ | $\alpha_3$ | $\alpha_4$ | $\alpha_5$ | $\alpha_6$ | FIC | AIC | BIC |
| 1 | 0.5 | - | - | - | - | - | 87 | 93 | 91 |
| 2 | 0.5 | 0.4 | - | - | - | - | 85 | 81 | 89 |
| 3 | 0.3 | 0.3 | 0.2 | - | - | - | 81 | 72 | 78 |
| 4 | 0.1 | 0.2 | 0.2 | 0.1 | - | - | 74 | 79 | 79 |
| 5 | 0.1 | 0.2 | 0.3 | 0.098 | 0.2 | - | 78 | 66 | 76 |
| 6 | 0.1 | 0.2 | 0.2 | 0.095 | 0.15 | 0.12 | 81 | 72 | 62 |
| $n = 2000$ | | | | | | | | | |
| True Order | Coefficients | | | | | | Scores | | |
| | $\alpha_1$ | $\alpha_2$ | $\alpha_3$ | $\alpha_4$ | $\alpha_5$ | $\alpha_6$ | FIC | AIC | BIC |
| 1 | 0.5 | - | - | - | - | - | 90 | 91 | 100 |
| 2 | 0.5 | 0.4 | - | - | - | - | 88 | 91 | 95 |
| 3 | 0.3 | 0.3 | 0.2 | - | - | - | 85 | 88 | 92 |
| 4 | 0.1 | 0.2 | 0.2 | 0.1 | - | - | 86 | 81 | 96 |
| 5 | 0.1 | 0.2 | 0.3 | 0.098 | 0.2 | - | 83 | 81 | 78 |
| 6 | 0.1 | 0.2 | 0.2 | 0.095 | 0.15 | 0.12 | 92 | 88 | 84 |
| $n = 3500$ | | | | | | | | | |
| True Order | Coefficients | | | | | | Scores | | |
| | $\alpha_1$ | $\alpha_2$ | $\alpha_3$ | $\alpha_4$ | $\alpha_5$ | $\alpha_6$ | FIC | AIC | BIC |
| 1 | 0.5 | - | - | - | - | - | 88 | 93 | 100 |
| 2 | 0.5 | 0.4 | - | - | - | - | 91 | 91 | 100 |
| 3 | 0.3 | 0.3 | 0.2 | - | - | - | 91 | 91 | 96 |
| 4 | 0.1 | 0.2 | 0.2 | 0.1 | - | - | 85 | 79 | 95 |
| 5 | 0.1 | 0.2 | 0.3 | 0.098 | 0.2 | - | 89 | 93 | 91 |
| 6 | 0.1 | 0.2 | 0.2 | 0.095 | 0.15 | 0.12 | 84 | 87 | 88 |

$\kappa = 2, \sigma = 1$

with that of the AIC and BIC.

We now report the simulations with focus on third central moment of errors. In this case, we assume $ALAR(p)$ to be a true model and treat corresponding $LAR(p)$ model as an alternative. Here, $p \in \{1, 2, ..., 7\}$. In each case, we recorded the number of correct selections made by AIC BIC and FIC. The results are reported in Table 2. It can be seen from Table 2 that the FIC outper-

Table 2. *Simulation results with focus on third central moment*

| $n = 500$ | | | | | | | | | | |
|---|---|---|---|---|---|---|---|---|---|---|
| Comparisons | Coefficients | | | | | | | Scores | | |
| | $\alpha_1$ | $\alpha_2$ | $\alpha_3$ | $\alpha_4$ | $\alpha_5$ | $\alpha_6$ | $\alpha_7$ | FIC | AIC | BIC |
| $LAR(1)$ vs $ALAR(1)$ | 0.5 | - | - | - | - | - | - | 91 | 83 | 98 |
| $LAR(2)$ vs $ALAR(2)$ | 0.5 | 0.3 | - | - | - | - | - | 85 | 72 | 97 |
| $LAR(3)$ vs $ALAR(3)$ | 0.3 | 0.3 | 0.3 | - | - | - | - | 94 | 78 | 98 |
| $LAR(4)$ vs $ALAR(4)$ | 0.1 | 0.3 | 0.1 | 0.2 | - | - | - | 96 | 77 | 98 |
| $LAR(5)$ vs $ALAR(5)$ | 0.1 | 0.2 | 0.3 | 0.098 | 0.18 | - | - | 90 | 66 | 95 |
| $LAR(6)$ vs $ALAR(6)$ | 0.1 | 0.2 | 0.2 | 0.098 | 0.095 | 0.1 | - | 86 | 65 | 95 |
| $LAR(7)$ vs $ALAR(7)$ | 0.1 | 0.2 | 0.2 | 0.098 | 0.095 | 0.15 | 0.12 | 91 | 68 | 96 |
| $n = 2000$ | | | | | | | | | | |
| Comparisons | Coefficients | | | | | | | Scores | | |
| | $\alpha_1$ | $\alpha_2$ | $\alpha_3$ | $\alpha_4$ | $\alpha_5$ | $\alpha_6$ | $\alpha_7$ | FIC | AIC | BIC |
| $LAR(1)$ vs $ALAR(1)$ | 0.5 | - | - | - | - | - | - | 96 | 79 | 98 |
| $LAR(2)$ vs $ALAR(2)$ | 0.5 | 0.3 | - | - | - | - | - | 91 | 78 | 98 |
| $LAR(3)$ vs $ALAR(3)$ | 0.3 | 0.3 | 0.3 | - | - | - | - | 93 | 81 | 100 |
| $LAR(4)$ vs $ALAR(4)$ | 0.1 | 0.3 | 0.1 | 0.2 | - | - | - | 98 | 79 | 100 |
| $LAR(5)$ vs $ALAR(5)$ | 0.1 | 0.2 | 0.3 | 0.098 | 0.18 | - | - | 94 | 79 | 99 |
| $LAR(6)$ vs $ALAR(6)$ | 0.1 | 0.2 | 0.2 | 0.098 | 0.095 | 0.1 | - | 91 | 77 | 97 |
| $LAR(7)$ vs $ALAR(7)$ | 0.1 | 0.2 | 0.2 | 0.098 | 0.095 | 0.15 | 0.12 | 91 | 78 | 100 |
| $n = 3500$ | | | | | | | | | | |
| Comparisons | Coefficients | | | | | | | Scores | | |
| | $\alpha_1$ | $\alpha_2$ | $\alpha_3$ | $\alpha_4$ | $\alpha_5$ | $\alpha_6$ | $\alpha_7$ | FIC | AIC | BIC |
| $LAR(1)$ vs $ALAR(1)$ | 0.5 | - | - | - | - | - | - | 92 | 85 | 100 |
| $LAR(2)$ vs $ALAR(2)$ | 0.5 | 0.3 | - | - | - | - | - | 91 | 75 | 99 |
| $LAR(3)$ vs $ALAR(3)$ | 0.3 | 0.3 | 0.3 | - | - | - | - | 95 | 82 | 99 |
| $LAR(4)$ vs $ALAR(4)$ | 0.1 | 0.3 | 0.1 | 0.2 | - | - | - | 93 | 81 | 99 |
| $LAR(5)$ vs $ALAR(5)$ | 0.1 | 0.2 | 0.3 | 0.098 | 0.18 | - | - | 94 | 75 | 100 |
| $LAR(6)$ vs $ALAR(6)$ | 0.1 | 0.2 | 0.2 | 0.098 | 0.095 | 0.1 | - | 92 | 69 | 98 |
| $LAR(7)$ vs $ALAR(7)$ | 0.1 | 0.2 | 0.2 | 0.098 | 0.095 | 0.15 | 0.12 | 97 | 74 | 100 |

$\kappa = 2, \sigma = 1$

forms AIC while it performs at par with the BIC. As the sample size increases, the performance of FIC as well as BIC improves though the improvement in the performance of BIC is slightly higher than that of FIC. Overall, the performance of FIC is comparable to that of BIC.

We repeated the above routine choosing various values of $\kappa$ below 1, above 1 and at 1 representing negatively skewed, positively skewed and symmetric Laplace autoregressive models respectively. Various sample sizes as well as model dimensions were chosen. The conclusions are consistent and coherent with those reported based on Tables 1 and 2. Hence, these simulations are not reported here to save the space.

### 7·2. *Method of moments* FIC

We demonstrate the focused selection between Gamma density with unknown scale $\alpha$ and shape $\beta$ and exponential density with unknown $\alpha$ and $\beta = 1$ using FIC based on method of moments. The moment function $h$ is taken to be $E\left(X, X^2\right)^{\mathrm{T}}$ where $X$ is a Gamma $(\alpha, \beta)$ random variable. The true model is chosen to be the exponential density with $\alpha = 2$. We generated 1000 samples of varying sizes. The proportions of correct selections made by method of moments FIC likelihood estimator based FIC AIC and BIC are reported in Table 3. The focus function was chosen to be the likelihood function of Gamma $(\alpha, \beta)$. It is seen that as sample size increases, FIC based on method of moments outperforms the other information criteria. Moreover, we do

Table 3. *Selection between Gamma and Exponential distributions*

| Sample size | FIC 1 | FIC 2 | AIC | BIC |
|---|---|---|---|---|
| 50 | 0.688 | 0.777 | 0.584 | 0.738 |
| 70 | 0.782 | 0.742 | 0.559 | 0.744 |
| 100 | 0.851 | 0.769 | 0.623 | 0.774 |
| 200 | 0.965 | 0.728 | 0.586 | 0.794 |
| 400 | 0.992 | 0.718 | 0.569 | 0.817 |
| 700 | 1 | 0.736 | 0.615 | 0.841 |
| 1000 | 1 | 0.752 | 0.564 | 0.843 |

FIC 1: FIC based on method of moments; FIC 2: FIC based on maximum likelihood estimators not observe any significant improvement in performance of other information criteria as sample size increases except in FIC based on method of moments. The possible explanation for this lies in the bias-variance trade-off in estimated focus function under local misspecification. Maximum likelihood estimators are highly sensitive to even moderate distributional misspecification. Moreover, in case of selection between narrow and wide models only without any other intermediate submodels, the maximum likelihood estimator of the mean squared error of focus function involves only the variance component under wide model and only the bias component under narrow model. See Chapter 6, page 147 of Claeskens and Hjort (2008). Thus, this FIC fails to take into consideration the effect of bias-variance trade-off. We report below the squared bias and variance of estimated focus function averaged over 1000 samples. A strikingly increasing trend

Table 4. *Bias-variance trade-off*

| Sample size | Bias 1 (narrow) | Bias 1 (wide) | Var 1 (narrow) | Var 1(wide) | Bias 2(narrow) | Var 2 (wide) |
|---|---|---|---|---|---|---|
| 50 | 1.1784 | 0.0070 | 0.1320 | 1.3309 | 0.044 | 0.058 |
| 70 | 1.9656 | 0.0092 | 0.2535 | 2.4152 | 0.0807 | 0.1055 |
| 100 | 3.2621 | 0.011 | 0.4916 | 4.3164 | 0.1380 | 0.1940 |
| 200 | 9.5997 | 0.01505 | 1.8412 | 14.6003 | 0.5793 | 0.7083 |
| 400 | 16.6489 | 0.02714 | 7.2311 | 30.9618 | 2.2587 | 2.6824 |
| 700 | 35.06933 | 0.4236 | 21.8620 | 61.2561 | 6.7109 | 7.9056 |
| 1000 | 96.4821 | 0.6639 | 44.4788 | 101.094 | 14.3888 | 16.0574 |

Bias 1 (narrow): Squared bias under narrow model based on method of moments; Bias 1 (wide): Squared bias under wide model based on method of moments; Var 1 (narrow): Variance under narrow model based on method of moments; Var 1 (wide): Variance under wide model based on method of moments; Bias 2 (narrow): Squared bias under narrow model based on maximum likelihood estimators; Var 2 (wide): Variance under wide model based on maximum likelihood estimators is seen in the bias and variance of moment estimator of the focus. The bias of moment estimators of focus in the narrow model is seen to be higher than that in the wide model and vice-versa for the variance of the estimated focus. In case of maximum likelihood based estimators of focus function, there is no role of bias in the wide model and of variance in the narrow model. Moreover, although the bias under narrow model and variance under wide model increase with sample size even in this case, the difference between them does not vary considerably when contrasted against the corresponding values for moment estimators. Hence, even for the data of size as large as 1000, maximum likelihood based FIC may not be able to distinguish correctly between narrow and wide models. For method of moment based FIC the dominance of bias in the narrow model is compensated by the dominance of variance in the wide model and vice-versa.

The above experiment was repeated with different values of $\alpha$ and with the focus on the moment generating function. The results were found to be favoring the method of moments based focused model selection.

## 8. CONCLUSIONS

We conclude that the FIC is a better, theoretically more sound alternative to classical likelihood based information criteria like AIC and BIC. The local misspecification setting is crucial in the focused model selection. Thus, we insist upon and recommend focused model selection by means of an FIC using non-likelihood, distribution-free and robust $M$-estimators which may utilize the bias-variance trade-off in the estimated focus functions while selecting the models. Hansen (2005) also recommends the focused model selection using such estimators. Our theoretical investigations and simulation results resonate with his recommendations.

## 9. DISCUSSIONS

### 9·1. *Frequentist model averaging using $M$-estimators*

The model selection is a crucial stage in statistical modeling and analysis. Quite often in practice, the final model chosen by an appropriate selector is treated as the true model as if known a priori and effect of model selection stage is ignored. In order to device valid and optimal post model selection inference, Hjort and Claeskens (2003) proposed the frequentist model averaging as an alternative to Bayesian model averaging of Hoeting et al. (1999). Instead of working with single best model, the model averaging framework extracts the "good" of all the candidate models with the resultant estimator as the random weighted average of the estimators based on each submodel. The model averaged estimator is shown to have smaller asymptotic mean squared error as opposed to the estimators based on the single model.

The model averaging setup of Hjort and Claeskens (2003) relies on maximum likelihood estimators under the assumption of independence of observations. By means of our asymptotic theory, it is possible to device the model averaging using any $M$-estimators in the dependent setup. The results of Hjort and Claeskens may then be viewed as a special case.

Let $\mathscr{M}$ denote the collection of all submodels $\mathcal{S}$. Let $w_n(\mathcal{S})$ denote random weight functions based on data of size $n$. The model averaged estimator, denoted by $\hat{\mu}_{\text{avg}}$ is given as

$$\hat{\mu}_{\text{avg}} = \sum_{\mathcal{S} \in \mathscr{M}} w_n(\mathcal{S})\, \hat{\mu}^{\mathcal{S}}_{t-1}.$$

The common choices of weights are the smoothed versions AIC, BIC and FIC as proposed by Hjort and Claeskens (2003). Assume that $\sum_{\mathcal{S} \in \mathscr{M}} w_n(\mathcal{S}) = 1$. In view of our Theorem 3, it can be shown on the lines of Hjort and Claeskens (2003) that

$$n^{1/2}\left(\hat{\mu}_{\text{avg}} - \mu^{\mathcal{S}}_{t-1}\right) \to \sum_{\mathcal{S} \in \mathscr{M}} w_n(\mathcal{S})\, N_1\left(m_{\mathcal{S}}^{\mathrm{T}}\eta, V_{\mathcal{S}}\right) \text{ in distribution.}$$

where $m_{\mathcal{S}}$ and $V_{\mathcal{S}}$ are given by (10) and (11) respectively. Using the above result, it is possible to construct optimum confidence intervals for $\mu_{\text{avg}}$ using any $M$-estimators, maximum likelihood estimators being a special case.

### 9·2. FIC *using stochastic vector valued focus functions*

For the sake of simplicity, we have assumed throughout the present discussions that conditional on $\mathscr{F}_{t-1}$, $\mu_{t-1}$ denotes a deterministic (non-stochastic) real-valued focus function.

However, the setup of functional calculus that we adopt enables immediate generalizations the above results for stochastic vector valued focus functions. For an $M$-functional $T$ in (1), let $\{\mu\left(T(F_n)\right)\}$ denote a sequence of focus functionals taking values in $\mathbb{R}^r, r \geq 1$. Assume that for every $n \geq 1$, $\mu\left(T(F_n)\right)$ is Hadamard differentiable at $T(F_\theta)$. The chain rule for Hadamard derivatives (see Van der Vaart 2000, Chapter 20, Theorem 20.9) then implies that the composite functional $\mu \circ T$ is Hadamard differentiable at $F_\theta$. Let $d_{F_0} = d_{F_{\theta_0}}\mu\left(T(F_n)\right)$ denote $k \times r$ Hadamard differential of $\mu$ at $F_{\theta_0}$. Applying functional delta method to (4) shows that $n^{1/2}\left(\mu\left(T(F_n)\right) - \mu\left(T(F_{\theta_0})\right)\right)$ converges in distribution under local misspecification to an $r$-variate Normal random vector with mean $-d_{F_0}^{\mathrm{T}} b(\delta)$ and dispersion matrix $d_{F_0}^{\mathrm{T}} D\Sigma D^{\mathrm{T}} d_{F_0}$. The FIC in this case is a consistent estimator of $||d_{F_0}^{\mathrm{T}} b(\delta)b(\delta)^{\mathrm{T}} d_{F_0}|| + ||d_{F_0}^{\mathrm{T}} D\Sigma D^{\mathrm{T}} d_{F_0}||$.

### 9·3. *Role of likelihood in model selection*

The asymptotic bias of $M$-estimators in our Theorem 2 involves the gradient of log-likelihood (score vector) although the estimating function is not necessarily specified in terms of likelihood. Does this fact contradict the very philosophy of robust, "distribution free" estimation? Why should a model selection criterion involve the likelihood at all if the underlying parametric estimation is not based on it? To this end, we reiterate that the parametric model selection here is perceived as the selection of probability models. It is possible in some cases to design model selection strategy in terms of the estimating functions alone without aids of any probabilistic specifications of data generating mechanisms. For instance, assuming the innovation probability distribution to be *correctly specified*, the variable selection in linear models or order selection in autoregressive models can be handled only by conditional expectation of response variable. In such cases, as advocated by DiTraglia (2016), the local misspecification can be introduced only in the underlying moment conditions. The FIC that he develops therefore does not involve any distributional attributes such as likelihood unlike the expressions that we obtain. However, not all the model selection problems can be handled in this way. If it is to agree upon that model selection is essentially about the choice of probability models, we believe that the misspecification should be introduced in probability measures instead of the moment conditions alone. Hence, it should not be very surprising or unjustified to let the likelihood enter the selection criterion. The reasons for working with $M$-estimators are associated with the robustness aspects of estimation. The local misspecification only partially serves the purpose of robustness since the parametric structural form of the data generating model needs to be specified at the outset. The data generating probability density is intractable in general, in which case, fully robust model selection procedures need to be developed. We provide the guidelines for the same below.

### 9·4. *Robustness aspects of* FIC

The motivation behind Huber's $M$-estimators is the qualitative robustness that may safeguard the inferential procedures against the gross errors in the data and moderate deviations from the assumed probability models. The $M$-estimators with bounded estimating function $\psi\left(.,\theta\right)$ have bounded gross-error sensitivity and hence are $B$ (bias)-robust. See Hampel et al. (2011), Chapter 4. However, bounding gross error sensitivity increases the asymptotic variance of estimators. Likelihood estimators, albeit asymptotically efficient, are non-robust in general. It is worth exploring how the trade-off between efficiency and robustness may be reflected in the focused selection of models. We demonstrated this phenomenon with moment estimators in Section 7.2. However, investigations of this nature call for further extensive theoretical as well as empirical analyses. The relative performance of various $M$-estimators in focused model selection in presence of outliers also needs to be evaluated and suitable estimators for model selection should be determined in this context. It is also possible to device model selection and averaging methods using optimal $B$-robust estimators of Hampel which achieve simultaneously the efficiency and the robustness. We leave these problems for future work.

The assumption of local misspecification partially serves the purpose of robust procedures, in a sense that it allows the deviations of the data generating model from the true, fixed model and yield bias-optimal estimators. However, the essence of "infinitesimal" robust procedures advocated by Hampel lies in *full neighbourhood systems* of Rieder (2012), Chapter 4 wherein the real distribution $P^*$ may be any member of full neighbourhood $\mathcal{U}_\theta$ of $P_\theta$ where $\mathcal{U}_\theta$ is a collection of all the open balls about $P_\theta$ in the space of all probability measures equipped with the appropriate topology. For any $k$-dimensional tangent $\zeta \in \mathscr{L}^2(P_\theta)$, the *simple perturbations* of $P_\theta$ along $\zeta$ are defined as

$$dP^*(\zeta, \eta) = \left(1 + n^{-1/2}\eta^{\mathrm{T}}\zeta\right) dP_\theta, \eta \in \mathbb{R}^k. \tag{16}$$

Relation (16) can be viewed as generalization the parametric local misspecification setting to a truly robust, distribution-free case. Lemma 4.2.4 and proposition 4.3.6 of Rieder (2012) are viewed as parallel versions of our Theorem 1 and 2 in the said infinitesimal robust setup. These results should facilitate the construction of robust FIC. The focus function in this case may be any real or vector valued Hadamard differentiable statistical functional defined on the space of probability measures. The functional delta method then yields the local asymptotics of an estimated focus. The true focus function is defined in terms of simple perturbations as $\mu_{\text{true}} = \mu(dP^*(\zeta, \eta))$. We discuss the theory and applications of robust FIC in our subsequent work.

## ACKNOWLEDGEMENT

The first author would like to thank the Council of Scientific and Industrial Research and University Grants Commission, Government of India for the award of Junior Research Fellowship.

## SUPPLEMENTARY MATERIAL

Supplementary material includes some auxiliary results, the proof of Theorem 1 and additional simulations using conditional least square estimators.

## APPENDIX

Suppose $F_n(s) = n^{-1}\sum_{i=1}^{n} I\{X_i \le s\}$ denotes the empirical distribution function of the data $X_1, X_2, ..., X_n$ while $F$ denotes the data generating distribution function. Let $B_F = \{B_F(s), s \in \mathbb{R}\}$ denote a zero mean Gaussian process (Brownian bridge) with covariance $E(B_F(t)B_F(s))$.

LEMMA A1. *Under A1-A4, the empirical process $Y_n = n^{1/2}(F_n - F)$ converges weakly in $D[-\infty, \infty]$ to $B_F$.*

*Proof.* The essential components of the proof are the convergence of the finite dimensional distribution of the underlying empirical process to that of the Brownian motion $B_F$, and the asymptotic tightness of the empirical process. See, for instance, Ossiander (1987), Theorem 3.1 and Goldie and Greenwood (1986 b).

*Convergence of finite dimensional distributions*: Assumptions A1-A4 show that the conditions of the central limit theorem (Lemma 2) of Deo (1973) are satisfied. Hence the finite dimensional distribution of $\{Y_n(s)\}$ converges to that of $B_F$.

*Asymptotic tightness*: The asymptotic tightness of the empirical process follows from Lemma 5.1 of Goldie and Greenwood (1986 b). In particular, our assumption A2 coincides with assumptions (i) and (ii) of Corollary 1.4 of Goldie and Greenwood. From the proof of the said corollary, it is clear that our A4 replaces assumption (iv) of Goldie and Greenwood.

Finally, it must be shown that the class $\mathscr{A} = \{I_{(-\infty,s]} : s \in [-\infty,\infty]\}$ of indicator functions on $(-\infty, s]$ satisfies the metric entropy condition (iii) of Corollary 1.4 of Goldie and Greenwood, or equivalently, the entropy integral condition of Theorem 3.1 of Ossiander (1987). Given $\epsilon > 0$, let $N_{[\,]}\left(\epsilon, \mathscr{A}, \mathscr{L}^2(F)\right)$ and $E_{[\,]}\left(\epsilon, \mathscr{A}, \mathscr{L}^2(F)\right) = \log N_{[\,]}\left(\epsilon, \mathscr{A}, \mathscr{L}^2(F)\right)$ denote respectively the bracketing number and the entropy of $\mathscr{A}$. See Van der Vaart (2000), Chapter 19, page 270 for the details. Following van der Vaart (2000), Ex. 19.6, it is seen that

$$N_{[\,]}\left(\epsilon^{1/2}, \mathscr{A}, \mathscr{L}^2(F)\right) \leq N_{[\,]}\left(\epsilon, \mathscr{A}, \mathscr{L}^2(F)\right) < 2/\epsilon. \tag{A1}$$

Hence,

$$\lim_{\epsilon \to 0} \sup \frac{\log E_{[\,]}\left(\epsilon, \mathscr{A}, \mathscr{L}^2(F)\right)}{\log(1/\epsilon)} < 1. \tag{A2}$$

Relation (A2) implies that condition (iii) of Corollary 1.4 of Goldie and Greenwood (1986 b) holds true for the class $\mathscr{A}$. Thus, the empirical process is asymptotically tight by virtue of the said lemma.

Since $\int_0^1 \log(2/\epsilon) d\epsilon < \infty$ , relation (A1) also implies that

$$\int_0^1 E_{[\,]}^{1/2}\left(\epsilon, \mathscr{A}, \mathscr{L}^2(F)\right) d\epsilon < \infty. \tag{A3}$$

Relation (A3) implies that the metric entropy condition of Theorem 3.1 of Ossiander (1987) is satisfied under our Assumptions A1-A4. The weak convergence of $Y_n$ thus now follows. $\square$

*Remark* A1. Billingsley (1968, Theorem 20.1) and Deo (1975) establish the asymptotic theory of empirical processes under a stronger condition of $\varphi$-mixing of stationary stochastic processes. The entropy integral conditions, on the other hand, enable the generalizations with weaker mixing conditions. Once the entropy integral conditions are established, the weak convergence of the tight empirical process is an immediate consequence of Donsker's theorem (van der Vaart, 2000, Theorem 19.5) leading to the shorter proofs than those by Billingsley and Deo.

Assume that a functional $T : D[-\infty, \infty] \to \mathbb{R}^k$ is Hadamard differentiable at $F$. Let $\mathbb{IC}(X_0) = \mathbb{IC}(X_0, T, F)$ denote the influence curve of $T$ at $F$. Assume that the influence curve is almost surely bounded. Let $V_0 = V\left(\mathbb{IC}(X_0)\right)$ denote $k \times k$ variance-covariance matrix of $\mathbb{IC}(X_0)$. Then,

THEOREM A1. *$n^{1/2}\left(T(F_n) - T(F)\right)$ converges (under $F$) in distribution to $k$-dimensional Normal random vector with mean vector $0_k$ and covariance matrix $V_0$.*

*Proof.* Since $T$ is Hadamard differentiable, the first order von Mises expansion of $T$ (Fernholz 1983, page 10, eq. 2.4) yields

$$n^{1/2}\left(T\left(F_n\right) - T(F)\right) = n^{1/2} dT_F\left(F_n - F\right) + n^{1/2} \text{Rem}\left(F_n - F\right), \tag{A4}$$

where $dT_F(Y)$ denotes $k \times 1$ Hadamard differential of $T$ with respect to $F$ at $Y$ and $\text{Rem}\left(F_n - F\right)$ denotes the remainder term. In view of proposition 4.3.3 of Fernholz, $\text{Rem}\left(F_n - F\right)$ can be seen to be a measurable element of $D[-\infty, \infty]$ although $F_n - F$ may not be. Since the empirical process $n^{1/2}\left(F_n - F\right)$ is shown to be asymptotically tight in the proof of Lemma A1,

$$n^{1/2} \text{Rem}\left(F_n - F\right) \to 0 \text{ in probability}. \tag{A5}$$

Moreover, the Hadamard derivative is again a measurable element of $D[-\infty, \infty]$. See Fernholz (1983), Lemma 4.4.1. Thus, $n^{1/2}\left(T\left(F_n\right) - T(F)\right)$ is measurable. Using relation (A5) in (A4) shows that the asymptotic distribution of $n^{1/2}\left(T\left(F_n\right) - T(F)\right)$ is same as that of $n^{1/2} dT_F\left(F_n - F\right)$. Since $T$ is Hadamard differentiable, functional delta method applied in conjunction with Lemma A1 shows that $n^{1/2}\left(T\left(F_n\right)-T(F)\right)$ converges in distribution to zero-mean $k$-dimensional Normal random vector. Further, Lemma 4.4.1 of Fernholz implies that

$$dT_F\left(F_n-F\right)=n^{-1}\sum_{t=1}^{n}\mathbb{IC}(X_t). \tag{A6}$$

The functional delta method implies that the asymptotic covariance matrix of $n^{1/2}dT_F\left(F_n-F\right)$ and in turn of $n^{1/2}\left(T\left(F_n\right)-T(F)\right)$ is same as that of $n^{-1/2}\sum_{t=1}^{n}\mathbb{IC}(X_t)$. Let $\mathbb{IC}_i(X_t)$ denote an $i$th element of $\mathbb{IC}(X_t)$ that represents the Hadamard differential of $T$ with respect to an $i$th component of $T(F)$, $i=1,2,...,k$. For any $i,j=1,2,....,k$, let $V_n(i,j)$ denote an $(i,j)$-th element of $V\left(n^{-1/2}\sum_{t=1}^{n}\mathbb{IC}(X_t)\right)$ while $V_0(i,j)$ denote an $(i,j)$-th element of $V\left(\mathbb{IC}(X_0)\right)$. To conclude the proof, we need to show that $V_n(i,j)\to V_0(i,j)$ for every $i,j=1,2,....,k$. Observe that

$$V_n(i,j)=n^{-1}\sum_{t=1}^{n}A_t(i,j)+n^{-1}\sum_{1\le t<s\le n}^{n}B_{t,s}(i,j)+n^{-1}\sum_{1\le s<t\le n}^{n}B_{s,t}(i,j), \tag{A7}$$

where

$$A_t(i,j)=E\left(\mathbb{IC}_i(X_t)\mathbb{IC}_j(X_t)\right)$$

and

$$B_{t,s}(i,j)=E\left(\mathbb{IC}_i(X_t)\mathbb{IC}_j(X_s)\right).$$

Since $\{X_n\}$ is a stationary process and influence curve is a measurable function of this process owing to relation (A6), $\{\mathbb{IC}_i(X_n)\}$ is also a stationary process for every $i=1,2,...,k$. Thus, $A_0(i,i)=n^{-1}\sum_{t=1}^{n}A_t(i,i)$ represents the variance of $\mathbb{IC}_i(X_0)$. Similarly, for every $i\neq j$, $n^{-1}\sum_{t=1}^{n}A_t(i,j)=A_0(i,j)$ represents the covariance between $\mathbb{IC}_i(X_0)$ and $\mathbb{IC}_j(X_0)$. Thus, we have for any $i,j=1,2....,k$,

$$A_0(i,j)=V_0(i,j). \tag{A8}$$

Without loss of generality, assume that $t<s$. Since the sequence $\{\mathbb{IC}_i(X_n)\}$ is stationary and strong mixing, Deo (1973), Lemma 3 implies that $\sum_{1\le t<s\le n}^{n}B_{t,s}(i,j)<\infty$ as $n\to\infty$ for every $i,j=1,2,...,k$. Thus,

$$n^{-1}\sum_{1\le t<s\le n}^{n}B_{t,s}(i,j)\to 0. \tag{A9}$$

Use A8 and A9 in A7 to conclude the proof. □